\documentclass[12pt,oneside,reqno]{amsart}

      \usepackage{amssymb}
      \usepackage{rotating}

\usepackage{amsfonts}
\usepackage{dsfont}
\usepackage{amsxtra}
\usepackage{pstricks}
\usepackage{pst-node}

\usepackage{index}
\newindex{A}{Aidx}{Aind}{Index of Symbols}
\newindex{B}{Bidx}{Bind}{Index of Notations}

      \theoremstyle{plain}
      \newtheorem{Theo}{Theorem}[section]
      \newtheorem{Lem}[Theo]{Lemma}
      \newtheorem{Prop}[Theo]{Proposition}
      \newtheorem{Cor}[Theo]{Corollary}
      
      \theoremstyle{definition}
      \newtheorem{Defi}[Theo]{Definition}
      \newtheorem{Exa}[Theo]{Example}

      \theoremstyle{remark}
      \newtheorem{Rem}[Theo]{Remark}

\def\ZZ{{\mathds Z}}
\def\NN{{\mathds N}}

\def\FF{{\mathds F}}

\def\hdp{\operatorname{Hdep}}
\def\dep{\operatorname{depth}}
\def\grade{\operatorname{grade}}
\def\pos{\operatorname{p}}

\def\lcm{\operatorname{lcm}}
\def\gcd{\operatorname{gcd}}

\def\sgr{\langle \alpha, \beta \rangle}

\def\lz{\vskip11pt}
\newcommand{\BE}{\hspace*{\fill} $\square$} 
      \makeatletter
      \def\@setcopyright{}
      \def\serieslogo@{}
      \makeatother
   
\catcode`\á=\active \def á{\'a}
 \catcode`\Á=\active \def Á{\'A}
  \catcode`\ó=\active \def ó{\'o}
  \catcode`\é=\active \def é{\'e}
  \catcode`\ú=\active \def ú{\'u}
  \catcode`\í=\active \def í{\'{\i}}
  \catcode`\ñ=\active \def ñ{\~n}
  \catcode`\Ñ=\active \def Ñ{\~N}
  \catcode`\¿=\active \def ¿{?`}
  \catcode`\º=\active \def º{$^{\underline{o}}$}
  \catcode`\ª=\active \def ª{$^{\underline{a}}$}
  \catcode`\¡=\active \def ¡{!`}
  \catcode`\â=\active \def â{\^{a}}
  \catcode`\ê=\active \def ê{\^{e}}
  \catcode`\î=\active \def î{\^{\i}}
  \catcode`\ô=\active \def ô{\^{o}}
  \catcode`\û=\active \def û{\^{u}}
  \catcode`\ç=\active \def ç{\c{c}}
  \catcode`\ü=\active \def ü{\"{u}}
  \catcode`\ö=\active \def ö{\"{o}}

 \newenvironment{dem}{\noindent{\it Proof}.}{\qed}

\newcommand{\mm}{\mathfrak{m}}

\begin{document}

%



   \author{Julio Jos\'e Moyano-Fern\'andez}
   \address{Institut f\"ur Mathematik, Universit\"at Osnabr\"uck. Albrechtstra\ss e 28a, D-49076 Osnabr\"uck, Germany}
   \email{jmoyano@math.uni-osnabrueck.de}


   \author{Jan Uliczka}

   \address{Institut f\"ur Mathematik, Universit\"at Osnabr\"uck. Albrechtstra\ss e 28a, D-49076 Osnabr\"uck, Germany}


   \email{juliczka@uos.de}
   

   \title[Hilbert depth of graded modules over polynomial rings]{Hilbert depth of graded modules over polynomial rings in two variables}


   \begin{abstract}
 In this article we mainly consider the positively $\ZZ$--graded polynomial ring $R=\FF[X,Y]$ over an 
arbitrary field $\FF$ and Hilbert series of finitely generated graded $R$--modules. The central result is
 an arithmetic criterion for such a series to be the Hilbert series of some $R$--module of positive depth. 
In the generic case, that is $\deg(X)$ and $\deg(Y)$ being coprime, this criterion can
be formulated in terms of the numerical semigroup generated by those degrees.    \end{abstract}



   \thanks{The first author was partially supported by the Spanish Government Ministerio de Educaci\'on y Ciencia (MEC)
grant MTM2007-64704 in cooperation with the European Union in the framework of the founds ``FEDER'', and by the Deutsche 
Forschungsgemeinschaft (DFG)}




   \maketitle



\section{Introduction and Review}

We want to investigate how some of the results of  \cite{ju} for Hilbert series of finitely generated graded 
modules over the standard $\ZZ$--graded polynomial ring can be generalised to the case where the ring of 
polynomials is endowed with an arbitrary positive $\ZZ$--grading.

Let  $R= \FF[X_1, \ldots , X_n]$ be the positively $\ZZ$--graded polynomial ring over some field $\FF$, i.~e.~each $X_i$ has degree $d_i \geq 1$ for every $i= 1, \ldots , n$. Moreover, let $M$ be a finitely generated graded 
$R$--module. Every homogeneous component of $M$ 
is a finite--dimensional $\FF$--vector space, and since $R$ is positively graded and $M$ is finitely generated, 
$M_j =0$ for $j \ll 0$. Hence the {\em  Hilbert function} of $M$
\[ H(M,-): \ZZ \to \ZZ,~ j \mapsto \dim_{\FF}\left(M_j \right), \] 
is a well--defined {\em  integer Laurent function} (see \cite{kr}, Definitions 5.1.1 and 5.1.12).
 The formal Laurent series {\em  associated to} $H(M,-)$
\[ H_M (t) = \sum_{j \in \ZZ} H(M,j) t^j 
=\sum_{j \in \ZZ} \left( \dim_{\FF} M_j \right) t^j \in  \ZZ[\![t]\!][t^{-1}] \]
is called the {\em  Hilbert series} of $M$. Obviously it has no negative coefficients; such a series will be
called {\em nonnegative} for short. 

By the Theorem of Hilbert--Serre (see \cite[Thm.~4.1.1]{ls}), $H_M$ may be written as a fraction of the form
\[ \frac{Q_M(t)}{\prod_{i=1}^n \left(1-t^{d_i} \right)},  \]
with some $Q_M \in \ZZ[t, t^{-1} ]$. As a consequence of this theorem and a well--known result 
in the theory of generating functions, see Proposition 4.4.1 of \cite{st}, there exists a quasi--polynomial 
$P$ of period $d:= \lcm(d_1 , \ldots ,d_n)$ such that $\dim_{\FF}\big( M_j \big) = P(j)$ for $j \gg 0$. 

The ring $R$ is $^*$local, that is, it has a unique maximal graded ideal, namely $\mm := (X_1 , \ldots, X_n)$. 
The {\em depth} of  $M$ is defined as the maximal length of an $M$--regular sequence in $\mm$, i.~e.~the grade of $\mm$ on $M$, 
and denoted by $\dep(M)$ rather than $\grade(\mm,M)$. This deviation from the standard terminology, where ``depth'' is used 
exclusively in the context of true local rings, may be justified by the fact that $\grade(\mm,M)$ agrees with 
$\dep\left(M_{\mm}\right)$, see \cite[Prop.~1.5.15]{bh}.

It is easy to see that (contrary to the Krull dimension) the depth of a module $M$ is not encoded in its
Hilbert series. Therefore it makes sense to introduce 
\[ \hdp(M) := \max \left\{ r \in \NN ~\left|~ \begin{array}{l} 
\mbox{There is a f.~g.~gr.~$R$--module $N$} \\ \mbox{with $H_N=H_M$ and $\dep(N) = r$.} \end{array}
 \right. \right\}; \]
this number is called the {\em Hilbert depth} of $M$. 

If the ring $R$ is standard graded, then $\hdp(M)$ turns out to coincide with the arithmetical invariant
 \[ \pos(M) := \max \left\{ r \in \NN ~|~ (1-t)^r H_M(t)~\mbox{is nonnegative} \right\}, \]
called the {\em positivity} of $M$, see Theorem 3.2 of \cite{ju}. The inequality  $\hdp(M) \leq \pos(M)$ follows from general 
results on Hilbert series and regular sequences. The converse can be deduced from the main result of \cite{ju}, Theorem 2.1, which 
states the existence of a representation
\[ H_M(t) = \sum_{j=0}^{\dim(M)} \frac{Q_j(t)}{(1-t)^j} ~~~~\mbox{with nonnegative}~~ Q_j \in \ZZ[t,t^{-1}].\]

We begin our investigation by establishing a similar decomposition theorem for Hilbert series of modules
over any positively $\ZZ$--graded polynomial ring. This result has some consequences for the Hilbert depth, but
it does not lead to an analogue of the equation $\hdp(M) = \pos(M)$ -- the occurence of different factors 
in the denominator of $H_M$ complicating matters. In section 3 we restrict our attention to polynomial rings in
two variables. For this special case we deduce an arithmetic characterisation of positive Hilbert depth. This criterion, 
surprinsingly related to the theory of {\em numerical semigroups}, is the main result of our paper.

\section{Preliminary results}

\subsection{A decomposition theorem}

Let $\FF$ be a field. We consider the polynomial ring $R= \FF[X_1, \ldots , X_n]$,  endowed with a general positive 
$\ZZ$--graded structure, i.~e.~$\deg(X_i) = d_i \geq 1$ for every $i= 1, \ldots , n$, and let $M$ be a finitely generated 
graded $R$--module. The Hilbert series of $M$ admits a decomposition analogous to that in the standard graded case 
(cf.~\cite[Thm.~2.1]{ju}). This can be proven using certain filtrations, similar to the argument in the proof of \cite[Prop.~2.13]{bku}: 
\begin{Theo} \label{2T1}
 Let $M$ be a finitely generated graded module over the positively graded ring of polynomials $R$. The Hilbert 
series of $M$ can be written in the form
\begin{equation} \label{2eq1}
 H_M(t) = \sum_{I \subseteq \{1, \ldots , n\}} \frac{Q_I(t)}{\prod_{i \in I} \left(1-t^{d_i} \right)} 
\end{equation}
with nonnegative  $Q_I \in \ZZ[t,t^{-1}]$.
\end{Theo}

\dem~ We use induction on $n$. For $n=0$, i.~e.~$R=\FF$, the module $M$ is just a finite--dimensional graded 
vector space, and hence $H_M$ itself is a Laurent polynomial. 
Assuming the claim to be true for modules over the polynomial ring in at most $n-1$ indeterminates, we consider a finitely 
generated module $M$ over $R = \FF[X_1, \ldots , X_n]$ and define a descending sequence of submodules
 $U_i$ of $M$ by $U_{n+1} :=M$ and 
\begin{equation} \label{2eq2}
 U_i:= \left\{ m \in U_{i+1} ~|~ X_{i}^jm = 0 ~ \mbox{for some}~ j>0 \right\}
\end{equation}
for $i= n , \ldots , 1$. Then for each $i$ we have a short exact sequence
\[ 0 \longrightarrow U_i \longrightarrow U_{i+1} \longrightarrow \underbrace{U_{i+1}/U_i}_{=:N_i} \longrightarrow 0 \]
and therefore
\[ H_{U_{i+1}} = H_{U_i} + H_{N_i}. \]
Combining these equations yields
\begin{equation} \label{2eq3}
 H_M = H_{U_1} + \sum_{i=1}^n H_{N_i}.
\end{equation}
Among all these Hilbert series, the first one is harmless, because it is easy to see that $U_1$ coincides with the local cohomology  
$H_{\mm}^0(M)$, and so it has finite length, see \cite[Prop.~3.5.4]{bh}. Therefore it is enough to show that each series $H_{N_i}$ 
admits a decomposition of the form (\ref{2eq1}). By construction, $X_i$ is not a zerodivisor on $N_i$, thus we have further exact 
sequences
\begin{equation*} 
0 \longrightarrow N_i(-d_i) \stackrel{ \cdot X_i}{\longrightarrow} N_i \longrightarrow 
\underbrace{N_i/X_iN_i}_{=:V_i} \longrightarrow 0,
\end{equation*}
and it follows that
\[ H_{V_i}(t) = \left(1-t^{d_i} \right) H_{N_i}(t). \]
Since $X_i$ acts trivially on $V_i$, we may regard $V_i$ as a module over the ring of polynomials 
$\FF[X_1, \ldots , \widehat{X_i}, \ldots , X_n]$ in $n-1$ indeterminates. By the induction hypothesis we can write the 
corresponding Hilbert series in the form
\[ H_{V_i}(t) = \sum_{I \subseteq \{1, \ldots , \hat{\imath}, \ldots , n\}} \frac{Q_I(t)}{\prod_{j \in I} \left(1-t^{d_j} \right)}. \]
Division by $1-t^{d_i}$ yields a presentation of the required form for the Hilbert series of 
$N_i$. \BE \medskip

A formal Laurent series admitting a decomposition of the form (\ref{2eq1}) will be called
\emph{$(d_1, \ldots, d_n)$--decomposable}.  Note that such a decomposition is
by no means unique. We define an important invariant of such a series:
\begin{equation}
\nu(H) := 
 \max \left\{r \in \NN ~\left|~
\begin{array}{l} \mbox{$H$ admits a decomp.~
  of form (\ref{2eq1}) } \\ \mbox{with}~ Q_I = 0 ~\mbox{for all}~ I ~\mbox{with}~ |I| <r.
\end{array} \right.\right\}
\end{equation}

It is easily seen that any $(d_1, \ldots, d_n)$--decomposable series $H$ is in fact the Hilbert 
series of some finitely generated $R$--module: Choose a decomposition of $H$ with $Q_I(t) =\sum_{k=p_I}^{q_I}h_{I,k}t^k$, and
 write $J_I$ for the ideal generated by the $X_i$ with $i \notin I$, then the $R$--module
\begin{equation} \label{2eq4}
 N:= \bigoplus_{I \subseteq \{1, \ldots , n\}} \left( \bigoplus_{k=p_I}^{q_I}
 \Big( \left(R/J_I \right) (-k) \Big)^{h_{I,k}} \right)
\end{equation}
has Hilbert series $H$. Hence we have shown the following.
\begin{Cor} \label{2C1}
A formal Laurent series  $H \in \ZZ  [\![ t ] \!][t^{-1}]$ is the Hilbert series of a finitely generated  graded 
module over the ring 
 $\FF[X_1, \ldots , X_n]$ with $\deg(X_i)=d_i$ if and only if it is $(d_1, \ldots, d_n)$--decomposable.
\end{Cor}

\begin{Rem} \label{2R1}
 Note that this result is not the complete analogue of Corollary 2.3 in \cite{ju}, since it remains open whether 
\emph{any} nonnegative series of the form $\frac{Q(t)}{\prod_{i=1}^n \left(1-t^{d_i}\right)}$ is
$(d_1, \ldots, d_n)$--decomposable.
\end{Rem}

\subsection{Hilbert Depth and Positivity}\label{subsect22}

Let $R= \FF[X_1, \ldots , X_n]$ be positively $\ZZ$--graded with $\deg(X_i) =d_i \geq 1$, and let $d := \lcm(d_1, \ldots , d_n)$. 
Moreover, let $M$ be a finitely generated graded $R$-module with Hilbert series $H_M$. As in \cite{ju} we define the 
{\em Hilbert depth} of $M$ by

\[ \hdp(M) := \max \left\{ r \in \NN ~\left|~ \begin{array}{l} 
\mbox{There exists a f.~g.~gr.~$R$--module $N$} \\ \mbox{with $H_N=H_M$ and $\dep(N) = r$.} \end{array} \right. \right\}. \]

On the other hand, the notion of positivity has to be adjusted to our new situation, since in general there is no element of 
degree $1$ and a fortiori no $M$--regular sequence consisting of such elements, hence one cannot
expect a relationship between $\pos(M)$ and $\hdp(M)$. Instead of $\pos(M)$ we consider 
\[ \pos_d(M) := \max \big\{ r \in \NN ~|~ \left(1-t^d\right)^r H_M(t)~\mbox{is nonnegative} \big\}. \]
This is an upper bound for $\hdp(M)$ for the same reason as $\pos(M)$ is in the standard graded case: Since extension of the
 ground field does not affect either the depth of a module or its Hilbert series, we may assume that $\FF$ is infinite, and in 
this case a maximal $M$--regular sequence can be composed of elements of degree $d$. (This can be seen by considering $\bigcup_{\mathfrak{p} \in \mathrm{Ass}(M)} \mathfrak{p} \subseteq \mm$ in degree $d$: Since a vector space over an infinite field cannot be the union of proper subspaces, equality means that $\mm$ agrees with some $\mathfrak{p}\in \mathrm{Ass}(M)$ in degree $d$, and this implies $\dep (M)=0$).

It remains to check whether we still have  equality. An inspection of the proof in the standard graded case shows that it is 
advisable to consider a third invariant, namely
\[ \nu(M) := \nu(H_M), \]
which is well--defined by Theorem \ref{2T1}. Note that in \cite{bku} \emph{this} number is called the Hilbert depth. 

For any decomposition of the form (\ref{2eq1}) the $R$--module $N$ given in (\ref{2eq4}) has 
$ \dep(N)= \min\big\{ \big. |I| ~\big|~ Q_I \neq 0 \big\}$. Therefore we have an inequality
\begin{equation*} \label{3eq(-1)}
 \nu(M) \leq \hdp(M) \leq \pos_d(M).
\end{equation*}
In the standard graded case, \cite[Thm.~2.1]{ju} also yields $\pos(M) \leq \nu(M)$ and hence the equality of 
all three numbers, as already mentioned above. This reasoning cannot be carried over to the general situation,
since the denominators in (\ref{2eq1}) are different from $1-t^{d}$.

\subsection{The case $\hdp(M) = 1$}

The method used in the proof of Theorem \ref{2T1} also yields $\hdp(M)=\nu(M)$ in the special case $\hdp(M)=1$.

\begin{Prop} \label{prop:hd_nu}
Let $M$ be a finitely generated graded $R$-module with
$\hdp(M) \geq 1$. Then we have $\nu(M) \geq 1$.
\end{Prop}

\dem~ We may assume $\dep(M) = \hdp(M) \geq 1$. Define submodules $U_i$ of $M$ and quotients 
$N_i$ and $V_i$ as in the proof of Theorem \ref{2T1}. By assumption we have $U_1 \cong  H_{\mm}^0(M) = 0$, 
see \cite[Prop.~3.5.4]{bh}, so we get the equation 
\[ H_M =  \sum_{i=1}^n H_{N_i}. \]
 The relation $H_{V_i}(t) = \big(1-t^{d_i} \big) H_{N_i}(t)$
implies $\nu(N_i) \geq 1$, and hence $\nu(M) \geq \min_i \left\{ \nu(N_i) \right\} \geq 1$.  \BE
\medskip

Thus, $\hdp(M)$ and $\nu(M)$ coincide at least in three special cases:

\begin{Cor} 
Let $M$ be a finitely generated graded module over the positively graded polynomial ring $\FF[X_1, \ldots, X_n]$.
If $\hdp(M) \leq 1$ or $\hdp(M) = n$, then $\nu(M)= \hdp(M)$.
\end{Cor}
\dem~ For $\hdp(M) \leq 1$, this follows from $\nu(M) \leq \hdp(M)$ and  Proposition \ref{prop:hd_nu}, and
if $\hdp(M) =n$, then $H_M$ is the Hilbert series of a free  module over $R$. \BE

\section{The case of the polynomial ring in two variables}

From now on, we only consider modules over the ring $\FF[X,Y]$ with $\alpha:= \deg(X)$ and $\beta:= \deg(Y)$
being coprime; we may assume $ \alpha < \beta$.
 
We will deduce an arithmetic characterisation of positive Hilbert depth, which leads to an analogue of 
the equation $\nu(M) = \hdp(M) = \pos(M)$ at least for the special case of $\hdp(M)=1$.

The necessary condition

\begin{equation} \label{3eq0}
 \pos_d(M) = \pos_{\alpha \beta}(M) >0  
\end{equation}

alone is not sufficient, as the following example shows.

\begin{Exa} \label{Exmp31}
Let $R=\FF[X,Y]$ with $\alpha=3$ and $\beta =5$ and consider the module
$M:= R \oplus \big(R/\mm \big)(-1) \oplus \big(R/\mm \big)(-2)$.
We have $\nu(M) = \hdp(M) = 0$ and $\pos_{15}(M) =1$.
\end{Exa}
\dem~ It is easily checked that $H(M,n+15) \geq H(M,n)$ holds for any $n \in \NN$. This implies $\pos_{15}(M) \geq 1$, and by 
$H(M,1) = H(M,16)$ we obtain  $\pos_{15}(M) =1$.

Now we show that any $R$--module with the same Hilbert series as $M$ has to have depth $0$. Assume on the contrary that there 
is a module $N$ with $H_N =H_M$ and $\dep(N) >0$. As mentioned earlier (Subsection \ref{subsect22}), the field $\FF$ may be assumed to be infinite, and so 
$R$ contains an $N$--regular element $z$ of degree $15$. It turns out that such an element cannot exist due to the fact that, by elementary linear algebra, any $\FF$--linear map $N_1 \oplus N_2 \to N_7 \oplus N_{11}$ has a nontrivial 
kernel: First the element $X^5$ cannot be $N$--regular, as one sees by considering the map
\[
N_1 \oplus N_2 \to N_7 \oplus N_{11}, \hspace{15pt}(n_1,n_2) \mapsto \left(X^5n_2,X^{10}n_1 \right). 
\]
We may therefore assume $z=\lambda X^5 + \mu Y^3$ with $\lambda, \mu \in \FF$, $\mu \neq 0$, and consider the map
$f: N_1 \oplus N_2 \to N_7 \oplus N_{11}$ with
\[ 
f(n_1,n_2) = \left(\lambda X^2n_1+ \mu Yn_2, \mu Y^2n_1 -\mu X^3n_2 \right). 
\]
There is also a non--zero element $(a,b) \in \ker f$, i.~e.
\begin{eqnarray*}
\lambda X^2a + \mu Yb &=& 0 \\
\mu Y^2a - \mu X^3b &=& 0.
\end{eqnarray*}
We multiply the first equation by $X^3$, the second by $Y$ and add both of them. This yields
\[ \left(\lambda X^5 + \mu Y^3 \right) a =0, \]
hence $a =0$, since $\lambda X^5 + \mu Y^3$ was assumed to be $N$--regular. But this implies that
$b \neq 0$ is annihilated by powers of $X$ and $Y$, hence $H_{\mm}^0(N) \neq 0$, contradicting $\dep(N) >0$. \BE

\medskip

Using similar arguments as above, one can deduce additional necessary conditions for positive Hilbert depth of a 
module $M$ over $\FF[X,Y]$ with $\deg(X)=3$, $\deg(Y)=5$. Let $H_M(t) = \sum_n h_nt^n$, then the coefficients have to satisfy 
\begin{eqnarray} 
h_n + h_{n+1} &\leq& h_{n+6} + h_{n+ 10}, \label{3eq1}  \\
h_n + h_{n+2} &\leq& h_{n+12} + h_{n+ 5}, \label{3eq2}  \\
h_n + h_{n+4} &\leq& h_{n+9} + h_{n+ 10},  \label{3eq2a}  \\
h_n + h_{n+7} &\leq& h_{n+12} + h_{n+ 10}  \label{3eq3}
\end{eqnarray}
for all $n \in \ZZ$. Our next example shows that these additional inequalities are still not sufficient 
to ensure positive Hilbert depth: Let
\[ M:=  \big(R/\mm \big)(-1) \oplus R/(Y) \oplus \big(R/(Y)\big)\,(-7) \oplus \big(R/(Y)\big)\,(-8), \]
then its Hilbert series
\[ H_M(t) = 1+t+t^3 + \sum_{n=6}^{\infty} t^n \]
satisfies conditions (\ref{3eq0})--(\ref{3eq3}), but $\nu(M) =0$. Assume on the contrary that it is possible to decompose
$H_M$ into summands of the form $\frac{ p_it^i}{1-t^3}$ and $\frac{ q_jt^j}{1-t^5}$ with $p_i,q_j \in \NN$.
Since $h_5=0$ such a decomposition requires the summand $\frac{1}{1-t^3}$, but then the remainder
\[ \tilde{H}(t) = \sum_{n=1}^{\infty} \tilde{h}_nt^n := H_M(t) - \frac{1}{1-t^3} = t +  t^7 + t^8 + t^{10}+ \cdots \]
cannot be decomposed further, because $\tilde{h}_1=1$ and $ \tilde{h}_4 = \tilde{h}_6=0$; note that $\tilde{H}$ does not 
satisfy (\ref{3eq3}) for $n=-6$. Hence $\nu(M)=0$ and, by Proposition \ref{prop:hd_nu}, also $\hdp(M) =0$.

A computation using {\tt Normaliz} (see \cite{no}) reveals two other necessary conditions for positive Hilbert depth, namely
\begin{eqnarray} 
 h_n + h_{n+1} + h_{n+2} &\leq& h_{n+4} + h_{n+5} + h_{n+6}, \label{3eq4} \\
h_n + h_{n+2} + h_{n+4} &\leq& h_{n+5} + h_{n+7} + h_{n+9} \label{3eq4a}
\end{eqnarray}
for all $n \in \ZZ$. The example above does not satisfy the first one for $n=-1$. 

One might observe that the constants in conditions (\ref{3eq0}) -- (\ref{3eq4a}) which are not multiples of 3 or 5 are the
numbers 1, 2, 4 and 7. These are the only positive integers not contained in $\langle 3,5 \rangle := 3\NN_0 + 5 \NN_0$, the 
so-called \emph{numerical semigroup} generated by 3 and 5. This is not a mere coincidence: The necessary and sufficient 
conditions for positive Hilbert depth to be developed in the sequel will turn out to be closely related to the theory of 
numerical semigroups, so it seems advisable to recall some basic facts of this theory here.

\medskip

\subsection{Numerical semigroups generated by two elements}

Let $S$ be a sub--semigroup of $\NN_0$ such that the greatest common divisor of all its elements is equal to $1$. Then the subset 
$\NN_0 \setminus S$ has only finitely many  elements, which are called the \emph{gaps} of $S$. 
Such a semigroup is said to be \emph{numerical}. The smallest element $c=c(S) \in S$ such that $n \in S$ for all 
$n \in \NN$ with $n \ge c$ is called the \emph{conductor} of $S$. The number of gaps is called the \emph{genus} 
of $S$ and is denoted by $g(S)$.

We are interested in numerical semigroups generated by two elements. Let $\alpha, \beta \in \NN$ with
 $\gcd(\alpha, \beta)=1$; we write $\langle \alpha, \beta \rangle := \NN_0 \alpha + \NN_0 \beta$
and denote the set of gaps of this semigroup by $L$.

\begin{Lem}[Cf.~\cite{rg}, Prop.~2.13] \label{3L1}
The semigroup $\langle \alpha, \beta \rangle$ generated by two positive integers $\alpha, \beta$ with $\gcd(\alpha, \beta)=1$ 
is numerical. Its conductor and genus are given by
\[
c=c\left(\langle \alpha, \beta \rangle\right) = (\alpha -1) (\beta - 1)~~\mbox{and}~~ 
g\left(\langle \alpha, \beta \rangle\right) = \frac{c}{2}.
\]
\end{Lem}

\begin{Lem}[\cite{rosales}, Lemma 1] \label{3L2}
Let $e \in \ZZ$. Then $e \notin  \langle \alpha, \beta \rangle $ if and only if there exist 
$k, \ell \in \NN$ such that $e=\alpha \beta - k \alpha - \ell \beta$.
\end{Lem}

\begin{Cor} \label{3C1}
 Let $k, \ell \in \NN$ such that $1 \leq k < \beta$ and $1 \leq \ell < \alpha$, then 
$| \alpha \beta -k \alpha - \ell \beta | \in L$.
\end{Cor}
\dem~ This follows immediately from the preceding lemma, since we have either
$| \alpha \beta -k \alpha - \ell \beta| = \alpha \beta - k \alpha - \ell \beta$
or 
\[ | \alpha \beta -k \alpha - \ell \beta| = k \alpha + \ell \beta - \alpha \beta
 = \alpha \beta - (\beta- k) \alpha - (\alpha - \ell) \beta. \]
\BE 

\begin{Cor} \label{3C2}
 Any integer $n>0$ has a unique presentation
\[ n = p \alpha \beta - a \alpha - b \beta \]
with integers $p>0$, $0 \leq a < \beta$ and $0 \leq b < \alpha$.
\end{Cor}
\dem~ Since the gaps of $\sgr$ are covered by Lemma \ref{3L2}, we have to show the existence of the presentation for
integers $n = k \alpha + \ell \beta$ with $k, \ell \geq 0$. Let $k= q \beta + r$, $\ell = q' \alpha + r'$ with $0 \leq r <\beta$
and $0 \leq r' < \alpha$, then
\[ n= k \alpha + \ell \beta = (q + q'+2) \alpha \beta - (\beta-r) \alpha - (\alpha - r') \beta. \]
The uniqueness follows easily since $\gcd(\alpha, \beta)=1$. Let
\[ p \alpha \beta - a \alpha - b \beta = n = p' \alpha \beta - a' \alpha - b' \beta, \]
then
\[ \big( (p-p') \beta - a + a' \big) \alpha = (b - b') \beta, \]
so $\alpha$ has to divide $|b-b'| < \alpha$ and hence $b=b'$. But this implies $|p-p'| \beta = |a - a'| < \beta$ and
therefore $a=a'$ and $p=p'$ as well. \BE \lz

The presentation mentioned above will be of particular importance for the gaps of $\sgr$. In the sequel
we will frequently use the notation
\begin{equation}\label{3eq4b}
 e = \alpha \beta - a(e)\cdot \alpha - b(e) \cdot \beta .
\end{equation}

Let $n$ be a nonzero element of $S$. The set
\[ \mathrm{Ap}(S,n):=\{h \in S \mid h-n \notin S\} \]
is called the Ap\'ery set of $n$ in $S$.

\begin{Lem} [\cite{rg}, Lemma 2.4]\label{Lem:apery}
Let $S$ be a numerical semigroup and let $n \in S \setminus \{0\}$. Then
\[ \mathrm{Ap}(S,n)=\{0=w(0), w(1), \ldots , w(n-1)\}, \]
where $w(i) := \min \{ x \in S~|~ x \equiv i \!\mod n\}$  for $0 \leq i \leq n-1$.
\end{Lem}

\begin{Lem}Let $n \in S \setminus \{0\}$, $i \in \{0, \ldots, n-1\}$. Let $x$ be an integer congruent 
to $w(i)$ modulo $n$. Then $x \in S$ if and only if $w(i) \le x$.
\end{Lem}

\dem~ By assumption we have $x \equiv w(i) \equiv i \mathrm{~mod~} n$. 
Hence $x \in S$ implies $x \le w(i)$ by the definition of $w(i)$.
Conversely, let $x \in \mathds{Z}$ be such that $x \geq w(i)$. This implies 
$x-w(i)=\lambda n \geq 0$, and so $x=\lambda n + w(i) \in S$.
\BE

\begin{Lem} [\cite{si}, Theorem 1]\label{Lem:si}
Let $S = \sgr$, then
\begin{eqnarray*}
 \mathrm{Ap}(S,\alpha) &=& \left\{s \beta ~|~ s= 0, \ldots , \alpha-1\right\}, \\
 \mathrm{Ap}(S,\beta) &=& \left\{r \alpha ~|~ r= 0, \ldots , \beta-1\right\}.
\end{eqnarray*}
\end{Lem}

As an immediate consequence of the preceding lemmas we get 
\begin{Cor} \label{CorAp}
 Let $0 < r < \beta$ (resp.~$0 < s < \alpha$), and let $x$ be an integer such that $0 < x < r \alpha$ and 
$x \equiv r \alpha \! \mod \beta$ (resp.~such that $0 < x < s \beta$ and $x \equiv s \beta \! \mod \alpha$). 
Then $x$ is a gap of $\sgr$.
\end{Cor}

\subsection{Arithmetical conditions for  $\hdp(M) >0$ }

Next we introduce the announced necessary and sufficient condition for positive Hilbert depth. We begin by reconsidering the 
special case $\alpha =3$, $\beta =5$ and the inequalities mentioned above:
\begin{eqnarray*}
h_n & \leq & h_{n+15}, \\
h_n + h_{n+1} &\leq& h_{n+6} + h_{n+ 10},   \\
h_n + h_{n+2} &\leq& h_{n+12} + h_{n+ 5},    \\
h_n + h_{n+4} &\leq& h_{n+9} + h_{n+ 10},   \\
h_n + h_{n+7} &\leq& h_{n+12} + h_{n+ 10}, \\
h_n + h_{n+1} + h_{n+2} &\leq& h_{n+5} + h_{n+6} + h_{n+7}, \\
h_n + h_{n+2} + h_{n+4} &\leq& h_{n+5} + h_{n+7} + h_{n+9}.
\end{eqnarray*}
We note some observations on the structure of these inequalities.

\begin{enumerate}
 \item For each index $i$ on the left there are indices $j, j'$ on the right such that $i \equiv j \!\mod 3$,
$i \equiv j' \!\mod 5$ and $i < j,j'$.
\item In each inequality the constants appearing on the left hand side are gaps of the semigroup $\langle 3,5 \rangle$, and the 
difference of any two of these gaps is also a gap.
\item The constants appearing on the right hand side are either gaps of $\langle 3,5 \rangle$ or multiples of 3 or 5. 
\end{enumerate}
This motivates the following definition.

\begin{Defi} 
Let $\alpha, \beta >0$ be coprime integers and let $L$ denote the set of gaps of $\langle \alpha,\beta \rangle$.
An \emph{$(\alpha,\beta)$--fundamental couple} $[I,J]$ consists of two 
integer sequences $I=(i_k)_{k=0}^m$ and $J=(j_k)_{k=0}^m$, such that
\begin{enumerate}
\item[(0)] $i_0=0$.
\item[(1)] $i_1, \ldots , i_m, j_1, \ldots , j_{m-1} \in L$ and $j_0,j_m \leq \alpha \beta$.
\item[]
\item[(2)] $\begin{array}{lllll} i_k \equiv j_k &\!\!\!\mod \alpha &~\mbox{and}~& i_k < j_k & ~~\mbox{for} ~~ k= 0, \ldots ,m;\\
        j_k \equiv i_{k+1} &\!\!\!\mod \beta &~\mbox{and}~& j_k > i_{k+1} & ~~\mbox{for} ~~ k= 0, \ldots ,m-1;\\
        j_m \equiv i_{0}  &\!\!\!\mod \beta &~\mbox{and}~& j_m \geq i_{0}. &
       \end{array}$
\item[]
\item[(3)] $|i_k - i_{\ell}| \in L ~~\mbox{for} ~~1 \leq k < \ell \leq m$.
\end{enumerate}
The number $|I|=m+1$ will sometimes be called the \emph{length} of the couple.
The set of $(\alpha,\beta)$--fundamental couples will be denoted by $\mathcal{F}_{\alpha,\beta}$.
\end{Defi}

Fundamental couples will be discussed in more detail in the next subsection; here we just note the simplest nontrivial
example: 
\begin{Rem}
 Let $e$ be a gap of $\sgr$ with $e = \alpha \beta - a \alpha - b \beta$, then
$\big[(0,e),((\beta -a) \alpha, (\alpha -b) \beta) \big]$ form an $(\alpha,\beta)$--fundamental couple.
\end{Rem}

\begin{Rem}
 The number of $(\alpha, \beta)$--fundamental couples grows surprisingly with increasing $\alpha$ and $\beta$. We
give some examples:
\[ \begin{array}{r|r|c}
    S= \langle \alpha, \beta \rangle & \left|\mathcal{F}_{\alpha,\beta}\right| & g(S) \\ \hline
\langle 4,5 \rangle & 14 & 6 \\ \langle 4,7 \rangle & 30 & 9 \\ 
\langle 6,11 \rangle& 728 & 25 \\ \langle 11,13 \rangle & 104\,006 & 60
   \end{array} \]
\end{Rem}

The main result of this paper is the following theorem:

\begin{Theo} \label{3T1}
Let $R=\FF[X,Y]$ be the polynomial ring in two variables such that $\deg (X)= \alpha$, $\deg (Y)=\beta$ with 
$\gcd(\alpha, \beta)=1$. Let $M$ be a finitely generated graded $R$--module.
Then $M$ has positive Hilbert depth if and only if its Hilbert series $\sum_n h_n t^{n}$ satisfies 
the condition \medskip

{\rm ($\star$)} \hfill $\sum_{i \in I} h_{i+n} \le \sum_{j \in J} h_{j+n}$ ~~ 
\mbox{~for all~} $n \in \ZZ,~ [I,J] \in \mathcal{F}_{\alpha,\beta}. $ \hspace*{\fill}
\end{Theo}

It is easily seen that condition ($\star$) is indeed necessary for positive Hilbert depth. First we note some elementary
remarks concerning ($\star$).

\begin{Lem} \label{3L3}
a) Let $H,H' \in \ZZ[\![t]\!][t^{-1}]$ be nonnegative Laurent series satisfying condition {\rm ($\star$)}, then the same
 holds for $H+H'$ and $t^iH$, where $i \in \ZZ$. \\
b) The series $\frac{1}{1-t^{\alpha}}$ and $\frac{1}{1-t^{\beta}}$ satisfy condition {\rm ($\star$)}. \\
c) Let $H(t)=\sum_{n=b}^{\infty} h_nt^n$ be a nonnegative Laurent series satisfying condition {\rm ($\star$)}, then the same 
holds for $\frac{1}{1-t^{\alpha}} H(t)$ and $\frac{1}{1-t^{\beta}} H(t)$.
\end{Lem}

\dem~~Assertion a) is obvious, and b) is also clear in view of the definition of an $(\alpha,\beta)$--fundamental couple.
 To prove c) it is, by symmetry, enough to consider 
\[ \frac{1}{1-t^{\beta}} \cdot H(t) = \left( \sum_{n=0}^{\infty} t^{n \beta} \right) \cdot
   \sum_{n=b}^{\infty} h_nt^n  = \sum_{n=b}^{\infty} \Big( \sum_{k \geq 0} h_{n - k \beta} \Big) t^n \]
with a \emph{finite} inner sum. If $H$ fulfills condition ($\star$), then for each $(\alpha,\beta)$--fundamental couple
$[I,J]$ the inequalities
\[
 \sum_{i \in I} h_{n- k \beta +i}  \leq  \sum_{j \in J} h_{n - k \beta + j } 
\]
hold for all $n \in \ZZ$ and all $k \in \NN_0$. Summing up yields
\[ \sum_{k \geq 0} \sum_{i \in I} h_{n- k \beta +i}  \leq  \sum_{k \geq 0} \sum_{j \in J} h_{n - k \beta + j }, \]
and hence
\[ \sum_{i \in I} \sum_{k \geq 0}  h_{n- k \beta +i}  \leq  \sum_{j \in J} \sum_{k \geq 0}  h_{n - k \beta + j }, \]
as desired. \BE

Let $M$ be a finitely generated graded $R$--module with $\hdp(M) >0$. By Proposition \ref{prop:hd_nu}, its
Hilbert series admits a decomposition of the form
\[ H_M(t) = \frac{Q_2(t)}{(1-t^{\alpha})(1-t^{\beta})} + \frac{Q_X(t)}{1-t^{\alpha}} + \frac{Q_Y(t)}{1-t^{\beta}} \]
with nonnegative  $Q_2,Q_X,Q_Y \in \ZZ[t,t^{-1}]$. These three summands in the decomposition above fulfill
($\star$) by the previous lemma, and so does their sum $H_M$. Hence we have proven:

\begin{Prop} \label{3P1}
Let $M$ be a finitely generated graded $R$--module with $\hdp(M) >0$. Then $H_M$ satisfies condition {\rm ($\star$)}.
\end{Prop}

The proof of the converse is more elaborate and requires several steps.

First we show that we may restrict our attention to Laurent series $\sum_n h_nt^n$ whose coefficients eventually become 
(periodically) constant. Since the Hilbert series of an $R$--module $M$ with $\dim(M)=1$ is of this form, in the sequel 
such a series will be called a \emph{series of dimension} $1$ for short. 

Let $M$ be a finitely generated graded $R$--module with Hilbert series satisfying condition ($\star$). By Theorem \ref{2T1}, 
we have a decomposition of the form
\[ H_M(t) = Q_0(t) + \frac{Q_X(t)}{1-t^{\alpha}} + \frac{Q_Y(t)}{1-t^{\beta}} + \frac{Q_2(t)}{(1-t^{\alpha})(1-t^{\beta})} \]
with nonnegative $Q_0,Q_X,Q_Y,Q_2 \in \ZZ[t,t^{-1}]$. The reduction to the case of a series of dimension 1 is based 
on the following observation: For any $r>0$ the series
\begin{eqnarray}
 &&H'(t) := H_M(t) - t^{r \alpha \beta} \frac{Q_2(t)}{(1-t^{\alpha})(1-t^{\beta})} \label{3Eq4}\\
& =&Q_0(t) + \frac{Q_X(t)}{1-t^{\alpha}} + \frac{Q_Y(t)}{1-t^{\beta}} 
 + Q_2(t) \cdot \frac{ \sum_{j=0}^{r \beta -1} t^{j \alpha}}{1-t^{\beta}} =:\sum_n h'_n t^n  \nonumber
\end{eqnarray}
is of dimension $1$; therefore we have to show that for some $r$ the series $H'$ also satisfies condition {\rm ($\star$)}. To this end, 
choose $r \in \NN$ large enough such that  $r \alpha \beta > \deg(Q_0) + \alpha \beta$. Then $h'_n = h_n$ for $n < \deg(Q_0) + \alpha \beta$,
and so all inequalities of {\rm ($\star$)} which are influenced by $Q_0$ are valid by assumption. The remaining inequalities of {\rm ($\star$)} are valid since 
for indices $n >\deg(Q_0)$ the coefficients of $H'$ agree with those of the series
\[
\frac{Q_X(t)}{1-t^{\alpha}} + \frac{Q_Y(t)}{1-t^{\beta}} 
 + Q_2(t) \cdot \frac{ \sum_{j=0}^{r \beta -1} t^{j \alpha}}{1-t^{\beta}},
\]
and the latter satisfies condition {\rm ($\star$)}  by Lemma \ref{3L3}. Hence it remains to prove the following assertion:
 
\begin{Prop} \label{3P2}
 Let $H(t) = \frac{Q(t)}{(1-t^{\alpha} ) (1-t^{\beta} )}$ be a nonnegative Laurent series of dimension $1$. If $H$ satisfies
 condition {\rm ($\star$)},  then the series is $(\alpha, \beta)$--decompos\-able with $\nu(H) = 1$.
\end{Prop}

An idea how to prove this proposition is not far to seek: Without loss of generality let $H(t) = \sum_{n \geq 0} h_nt^n$ with 
$h_0 >0$. If we can show that at least one of the series
\[ H_{\alpha}(t) := H(t) - \frac{1}{1-t^{\alpha}} \,\,\,\,{\mbox{and}} \,\,\,\, H_{\beta}(t) := H(t) - \frac{1}{1-t^{\beta}}\]
is nonnegative and satisfies condition ($\star$), then a simple inductive argument would complete the proof.  

Indeed condition ($\star$) ensures that one of these series must be nonnegative:

\begin{Prop} \label{3P3}
 Let $\sum_{n=0}^{\infty}h_nt^n = \frac{Q(t)}{(1-t^{\alpha} ) (1-t^{\beta} )}$ be a nonnegative Laurent series 
 satisfying condition {\rm ($\star$)}. If $h_0 >0$, then at least one of the numbers
$ c_{\alpha} := \min \{ h_{r \alpha} ~|~ r>0 \}$ and $c_{\beta} := \min \{ h_{r \beta} ~|~ r>0 \}$
is positive.
\end{Prop}

\dem~ Let $c_{\alpha} = h_{k \alpha}$ and $c_{\beta} = h_{\ell \beta}$. Since $h_n \leq h_{n+ \alpha \beta}$ for all $n \in \ZZ$
and $h_0>0$ we may assume $0 < k < \beta$ and $0 < \ell < \alpha$. We consider 
\[ e := \alpha \beta - (\beta- k) \alpha - (\alpha-\ell) \beta. \] 
If this number is a gap of $\sgr$, then $\big( (0,e),(k \alpha, \ell \beta) \big)$ is a fundamental couple.
Hence we have
\[ 0 < h_0 \leq h_0 + h_e \leq h_{k \alpha} + h_{\ell \beta}, \]
and therefore at least one term on the right must be positive. Otherwise $e <0$ and, by Corollary \ref{3C1}, 
$-e= \alpha \beta - k \alpha - \ell \beta$ is a gap. In this case $\big( (0,-e),((\beta-k) \alpha, (\alpha -\ell) \beta) \big)$ is a 
fundamental couple, and hence the inequality 
\[  h_{n} + h_{n-e} \leq h_{n+(\beta -k) \alpha} + h_{n+(\alpha -\ell) \beta} \]
holds for all $n \in \ZZ$. For $n = e$ this yields
\[ 0 < h_0 \leq h_e + h_0 \leq h_{\ell \beta} +  h_{k \alpha}, \]
so again the right hand side must be positive. \BE \lz

The question whether condition ($\star$) is still valid for $H_{\alpha}$ or $H_{\beta}$
is more delicate. Subtraction of, say, $\frac{1}{1-t^{\alpha}}$ diminishes
all coefficients $h_{r \alpha}$ with $r \geq 0$ by 1; therefore all inequalities of ($\star$) containing such a coefficient
either on both sides or not at all are preserved. But there are (finitely many) inequalities where an index $r \alpha \geq 0$
on the right has a counterpart $r' \alpha <0$ on the left. We introduce a name for such an inequality.

\begin{Defi} \label{3D1}
 Let $H(t) = \sum_{n=0}^{\infty} h_nt^n$ be a formal Laurent series, and let $m \in \ZZ$ such that there exist $i \in I$, $j \in J$
 with $m + i <0$, $m+j \geq 0$ and $m+i \equiv m+j \equiv 0 \mod \alpha$. Then the inequality
\[\sum_{i \in I} h_{m+i} \leq \sum_{j \in J} h_{m+j}\]
is called \emph{$\alpha$--critical}. A \emph{$\beta$--critical} inequality is defined analogously.
\end{Defi}

Subtraction of $\frac{1}{1-t^{\alpha}}$ diminishes only the right hand side of an $\alpha$--critical inequality,
hence condition ($\star$) remains valid for $H_{\alpha}$ if and only if for $H$ all $\alpha$--critical inequalities are
strict, or likewise with $\beta$. Therefore we have to investigate whether for every Laurent series satisfying ($\star$)
the critical inequalities of at least one type hold strictly. This requires some technical machinery to be developed in the next
subsection.

\subsection{Fundamental and balanced couples}

The structure of the fundamental couples is widely determined by the following facts:
\begin{Lem} \label{3L4}
Let $i_1,i_2$ be gaps of $\sgr$, and let $a_k = a(i_k)$, $b_k=b(i_k)$ for $k=1,2$ denote the coefficients in the presentation of $i_k$ according to Corollary \ref{3C2}. Then

a) The difference $|i_1 - i_2|$ is a gap if and only if $(a_2-a_1)(b_2-b_1) <0$.

b) There exists a gap $j \geq i_1,i_2$ with $i_1 \equiv j \mod \alpha$ and $i_2 \equiv j \mod \beta$ if and only if
$a_1 \geq a_2 $ and $b_1 \leq b_2$, and this gap $j$ is uniquely determined to be 
$ \alpha \beta - a_2 \alpha - b_1 \beta$.

\end{Lem}

\dem~ a) We may assume $|i_1 - i_2| = i_1 -i_2 = (a_2-a_1) \alpha + (b_2-b_1) \beta$. 
If this number is a gap, then $a:=a_2 -a_1$ and $b:= b_2 -b_1$ must  bear different signs. For the converse note
that $|a| < \beta$ and $|b| < \alpha$, hence by Lemma \ref{3L2}
\[ i_1-i_2 = \alpha \beta - (\beta - a) \alpha + b \beta = \alpha \beta + a \alpha - (\alpha -b) \beta\]
is a gap if $a>0$ and $b<0$ or vice versa. 

b) Let $j:= \alpha \beta - a_2 \alpha - b_1 \beta$, then every solution of the congruence system in question is 
of the form $j + r \alpha \beta$ with some $r \in \ZZ$ (Chinese remainder theorem). By Lemma \ref{3L2}, $j$ is the 
only solution which is possibly a gap. On the other hand $j \geq i_1, i_2$ if and only if $a_1 \geq a_2 $ and $b_1 \leq b_2$, 
and in this case $j$ is indeed a gap.
\BE \medskip

The condition in part b) of the lemma will occur frequently in the sequel. Therefore we introduce a 
relation $\preceq$ as follows:
\begin{Defi}
 For gaps $i_1,i_2$ of the semigroup $\langle \alpha, \beta \rangle$ we define
\[ i_1 \preceq i_2 ~~ :\iff ~~ a(i_1) \geq a(i_2) ~ \wedge ~ b(i_1) \leq b(i_2) \]
and 
\[ i_1 \prec i_2 ~~ :\iff ~~ a(i_1) > a(i_2) ~ \wedge ~ b(i_1) < b(i_2). \]
\end{Defi}

Note that, deviating from the usual convention, $i_1 \prec i_2$ is a stronger assertion than just
$i_1 \preceq i_2 \,\wedge\, i_1 \neq i_2$.

Obviously $\preceq$ defines a partial ordering on the set of gaps. Together with the second part of Lemma \ref{3L4}  
this yields the announced structural result for fundamental couples.

\begin{Cor}  \label{3C3}
a) Let  $\left[(i_k)_{k=0}^m,(j_k)_{k=0}^m \right]$ be a fundamental couple, then 
$i_k \prec i_{k+1}$ for $k= 1, \ldots , m-1$. 

b) An $(\alpha,\beta)$--fundamental couple has length at most $\alpha$.

c) Let $i_1 \prec  i_2 \prec  \cdots \prec  i_m$ be gaps of $\sgr$, then there exists a unique sequence $J=(j_k)_{k=0}^m$ such
that $\big[(i_0:= 0,i_1, \ldots , i_m), (j_0, \ldots , j_m) \big]$ is a fundamental couple.

d) Let $L'=\{\ell_1, \ldots , \ell_m\}$  be a subset of $L:= \NN \setminus \sgr$ with $|\ell_n - \ell_p| \in L$ for
all  $n \neq p$. Then there exists a unique fundamental couple $[I,J]=\left[(i_k)_{k=0}^m,(j_k)_{k=0}^m \right]$ 
such that $L' = \{i_k~|~ 1 \leq k \leq m\}$.
\end{Cor}

\dem~ Part a) is  clear by Lemma \ref{3L4}, and so is c), since there are unique gaps $j_k$, 
$0<k<m$ such that $i_k \equiv j_k \!\mod \alpha$, $i_{k+1} \equiv j_k \!\mod \beta$ and $j_k > i_k, i_{k+1}$,
and the fundamental couple can only be completed by setting $j_0 :=  (\beta - a(i_1)  ) \alpha$ and 
$j_m :=  (\alpha - b(i_m)) \beta$. Part b) follows, because the integers $0< b(i_k) <  \alpha$, $k= 1, \ldots , m-1$ must be 
distinct. For d) one notes that all elements of $L'$ are $\prec$--comparable. Hence they can be ordered 
by this relation; it remains to apply the previous parts a) and c).   \BE \medskip

The fact that the series $\sum_{n\geq 0} t^{n \alpha}$ and $\sum_{n\geq 0} t^{n \beta}$ satisfy condition ($\star$) only depends
on the second requirement in the definition of a fundamental couple. This suggests this property to be the most important
for our purpose. We introduce a notion for couples of integer sequences with just this property:
\begin{Defi} 
Let $\alpha, \beta >0$ be coprime integers. 

a) An \emph{$(\alpha,\beta)$--balanced couple} $[I,J]$ consists of two 
integer sequences $I=(i_k)_{k=0}^m$ and $J=(j_k)_{k=0}^m$, such that
\[ \begin{array}{llll} i_k \equiv j_k \!\mod \alpha &~\mbox{and}~& i_k \leq j_k & ~~\mbox{for} ~~ k= 0, \ldots ,m;\\
        j_k \equiv i_{k+1}  \!\mod \beta &~\mbox{and}~& j_k \geq i_{k+1} & ~~\mbox{for} ~~ k= 0, \ldots ,m-1; \\
        j_m \equiv i_{0}  \!\mod \beta &~\mbox{and}~& j_m \geq i_{0}. & 
       \end{array}\]
The number $m+1$ will again be called the \emph{length} of the couple.

b) An  $(\alpha,\beta)$--balanced couple is called \emph{reduced}, if it satisfies the additional condition
\begin{eqnarray*}
 \min \left\{ j_{k-1} - i_k, j_k -i_k \right\} &< & \alpha \beta \hspace{1cm} \mbox{for}~ k = 1, \ldots , m \\
  \min \left\{ j_m - i_0, j_0 -i_0 \right\} &< & \alpha \beta  \\
 \min \left\{ j_k - i_k, j_k -i_{k+1} \right\} &< & \alpha \beta \hspace{1cm} \mbox{for}~ k = 0, \ldots , m-1 \\
\min \left\{ j_m - i_m, j_m -i_0 \right\} &<& \alpha \beta
\end{eqnarray*}
and the inequalities in a) hold strictly.
\end{Defi}

\begin{Rem}
a) Note that the length of a balanced couple is \emph{not} bounded
above, because there is no restriction against repetition of residue classes; even the same integer may appear
several times. 

b) Any fundamental couple except $[(0),(\alpha \beta) ]$ is also a reduced balanced couple.

c) Reducedness of the balanced couple $[I,J]$ implies that each $j_k$ is the smallest solution of the congruence system 
$x \equiv i_k \!\mod \alpha ~\wedge~ x \equiv i_{k+1} \!\mod \beta$ fulfilling the additional requirement
$x \geq i_k,i_{k+1}$. In particular, for any reduced balanced couple $[I,J]$, the sequence $J$ is uniquely determined by $I$.

\end{Rem}

\begin{Lem}  \label{3L6}
 Let $\left[(i_k)_{k=0}^m,(j_k)_{k=0}^m \right]$ be a reduced $(\alpha,\beta)$--balanced couple with nonnegative elements, $i_0=0$
and $m \geq 1$. 

a) At least one of the elements $i_1$, $i_m$ is a gap of $\sgr$. 

b) If $i_k,i_{k+1}$ are gaps with $i_k \preceq i_{k+1}$, then $i_k \prec i_{k+1}$.

c) If $i_1, \ldots , i_m$ are gaps such that $i_1 \preceq i_2 \preceq \cdots \preceq i_m$, then the couple
is fundamental.

d) Any reduced $(\alpha,\beta)$--balanced couple $\left[(0,i_1),(j_0,j_1)\right]$ is fundamental.
\end{Lem}

\dem~ a) The couple $[I,J]$ is reduced, therefore we have
\[ \min \left\{ j_m - i_0, j_0 -i_0 \right\} = \min \left\{j_m = s \beta, j_0 = r \alpha \right\} < \alpha \beta. \]
Since $i_1 \equiv j_0 \!\mod \beta$ and $i_m \equiv j_m \!\mod \alpha$ the claim follows immediately from Corollary 
\ref{CorAp}. 

b) If $a(i_k) = a(i_{k+1})$, then $i_k \equiv i_{k+1} \!\mod \beta$ and hence $j_k = i_k$, a contradiction;
the equality $b(i_k) = b(i_{k+1})$ is treated analogously.

c) By the previous part we even have $i_1 \prec  i_2 \prec  \cdots \prec  i_m$. Hence there exists a fundamental
couple with this $I$ by Corollary \ref{3C3}, and it is the only reduced balanced couple with this $I$.

d) This follows immediately from a) and c).
\BE

\begin{Lem}  \label{3L7}
 Let $\left[(i_k)_{k=0}^m,(j_k)_{k=0}^m \right]$ be an $(\alpha,\beta)$--balanced couple with
nonnegative elements, and let $i_k = \alpha \beta - a_k \alpha - b_k \beta$ be a gap of $\sgr$ for some $0 < k < m$. 
Then $(\alpha - b_k) \beta < j_k$ unless $i_{k+1}$ is a gap with $i_k \preceq i_{k+1}$, and 
vice versa~$(\beta - a_k) \alpha < j_{k-1}$ unless $i_{k-1}$ is a gap with $i_{k-1} \preceq i_k$.
\end{Lem}

\dem~ Let $i_{k+1} = p \alpha \beta - a_{k+1} \alpha - b_{k+1} \beta$ according to Corollary \ref{3C2}, then 
\[ j_k =  r \alpha \beta - a_{k+1} \alpha - b_{k} \beta \]
for some $ r \in \mathds{N}$. If $i_{k+1}$ is a gap then, by assumption, $i_k \not\preceq i_{k+1}$, and so,
as already mentioned in  the proof of part b) of Lemma \ref{3L4},
\[ j_k \geq 2 \alpha \beta - a_{k+1} \alpha - b_{k} \beta > \alpha \beta - b_k \beta; \]
otherwise $p>1$, and since
\[ (p-1) \alpha \beta - a_{k+1} \alpha - b_{k} \beta = p \alpha \beta - a_{k+1} \alpha - (\alpha +b_{k})\beta < i_{k+1}, \]
we have $r \ge p$ and therefore
\[ j_k \geq p \alpha \beta - a_{k+1} \alpha - b_{k} \beta \geq  2 \alpha \beta - a_{k+1} \alpha - b_{k} \beta 
> \alpha \beta - b_k \beta.\]
The second assertion can be proven analogously. \BE

The next result provides the key for showing Proposition \ref{3P2}. Its intricate proof is the technically most
challenging step on the way to our main result.

\begin{Theo} \label{Th:bl}
Let $H(t) = \sum_n h_nt^n$ be a nonnegative formal Laurent series satisfying condition {\rm($\star$)}. 
Then the inequality
\[ \sum_{i \in I} h_i \leq \sum_{j \in J} h_j \]
holds for any $(\alpha,\beta)$--balanced couple $[I,J]$.
\end{Theo}

\dem~ We may assume that $[I,J]$ is reduced: A perhaps necessary replacement of an $i \in I$ with $i + \alpha \beta$  or 
a $j \in J$ with $j - \alpha \beta$ is harmless since $h_n \leq h_{n + \alpha \beta}$ for all $n \in \ZZ$, while any elements 
$i_k = j_k$ or $j_k=i_{k+1}$ can be removed from $I$ and $J$ without affecting the inequality in question. Therefore we may
in particular assume $i_k \neq i_{k+1}$ for $k= 0, \ldots , m-1$.
Since $\left[ (i_k-x)_{k=0}^m, (j_k-x)_{k=0}^m \right]$ is a reduced balanced couple as well for any $x \in \ZZ$, 
we may also assume $\min I = 0$. Finally, we may shift the numbering of the elements in $I$ and $J$ such that $i_0 =0$.
Throughout this proof $a_k$ and $b_k$ denote the coefficients of $\alpha$ resp.~$\beta$ in the presentation of $i_k$ 
according to Corollary \ref{3C2}.

The proof uses induction on $m$, the case $m=0$ being trivial, while $m=1$ is covered by Lemma \ref{3L6}.

Let therefore $m \geq 2$ and assume that the result is already proven for balanced couples of length $\le m$. The general 
idea is to insert an auxiliary element $x$ into $I$ and $J$, which allows to split the amended couple into smaller balanced 
couples $[I',J']$ and $[I'',J'']$ with $I' \cup I'' = I \cup \{x\}$ and  $J' \cup J'' = J \cup \{x\}$. The inequalities
\[ \sum_{i \in I'} h_i \leq \sum_{j \in J'} h_j ~~\mbox{and}~~  \sum_{i \in I''} h_i \leq \sum_{j \in J''} h_j \]
then hold by the induction hypothesis, so we get our desired inequality by adding them and cancelling  $h_x$.

Since $[I,J]$ is reduced, at least one of the elements $j_0 = r \alpha$, $j_m = s \beta$ has to be less than $\alpha \beta$.
We distinguish three cases:

I) $j_0 <  \alpha \beta$ and $j_m \geq \alpha \beta$: In this case $i_1$ is a gap by Corollary \ref{CorAp}, while $i_m$ is not.
Let $M$ be the largest index $k$ such that $i_1, \ldots , i_k$ are gaps with $i_1 \prec i_2 \cdots \prec i_k$. Then $i_{M+1}$
is not a gap with $i_M \preceq i_{M+1}$. Hence Lemma \ref{3L7} implies that $x:= (\alpha - b_M) \beta < j_M$, 
and of course $x < \alpha \beta \leq j_m$ as well. Since $x \equiv i_M \equiv j_M \!\mod \alpha$ and $x \equiv j_m \!\mod  \beta$ 
we have two balanced couples
\begin{eqnarray} \label{3Eq5}
& &  \big[ (i_0, \ldots , i_M), (j_0, \ldots , j_{M-1}, x)\big] \\
&~~\mbox{and}~~&  \big[ (x, i_{M+1}, \ldots , i_m), (j_M,  \ldots , j_m)\big].\nonumber
\end{eqnarray}
Of these, the first one is already fundamental (by Lemma \ref{3L6}), while the second has length $m-M+1 \le m$, so the 
induction hypothesis can be applied. 

II) $j_m <  \alpha \beta$ and $j_0 \geq \alpha \beta$: This case is mirror--imaged to the first. Now $i_m$ is a gap 
and $i_1$ is not, so there is a smallest index $N$ such that $i_N, \ldots , i_m$ are gaps with $i_N \prec \cdots \prec i_m$.
Then $i_{N-1}$ is not a gap with $i_{N-1} \preceq i_N$, so $x:= (\beta - a_N) \alpha < j_{N-1}$ by Lemma \ref{3L7}, 
and also $x < \alpha \beta \leq j_0$. Since $x \equiv i_N \equiv j_{N-1} \!\mod \beta$ and $x \equiv j_0 \!\mod  \alpha$ we have
two balanced couples
\begin{eqnarray} \label{3Eq6}
& & \big[ (x, i_1, \ldots , i_{N-1}), (j_0, \ldots , j_{N-1})\big] \\
&~~\mbox{and}~~& \big[ (i_N, \ldots , i_m, i_0), (j_N,  \ldots , j_m, x)\big], \nonumber
\end{eqnarray}
the first being of length $N \le m$, and the second being a cyclic permutation of a fundamental couple. 

III) $j_0, j_m <  \alpha \beta$: In this case both $i_1$ and $i_m$ are gaps. We choose $M$ and $N$ as in case I) resp.~case II).
 If $M=m$, then the couple is fundamental and we are done. We may therefore assume $M <m$ and thus $N>M$. Two subcases can be
treated analogously to the cases above:

If $b_M  \geq  b_m$, then $x:=(\alpha-b_M)\beta \leq (\alpha-b_m)\beta=j_m$, so we may adopt the reasoning of
case I) and split $[I,J]$ into the couples given in (\ref{3Eq5}).

If $  a_N \geq   a_1$ then $x:=(\beta-a_N)\alpha \leq (\beta-a_1)\alpha=j_0$, so we may adopt the reasoning of
case II) and split $[I,J]$ into the couples given in (\ref{3Eq6}).

We may therefore assume that $b_M < b_m$ and $a_N < a_1$ and hence $a_m < a_1$ and $b_m > b_1$.

This case is treated recursively: Starting with
$p_0:= M$, $q_0:= N$, $u_1 := 1$ and $v_1 := m$ we construct two nonincreasing integer sequences $(p_r)$, $(v_r)$ 
and two nondecreasing integer sequences $(q_r)$, $(u_r)$ such that 
\begin{eqnarray}
 b_{p_{r-1}} < b_{v_r}~~ &\mbox{and} & ~~  a_{q_{r-1}} < a_{u_r}, \label{3Eqb} \\
p_{r-1} \geq u_r~~ &\mbox{and} & ~~  q_{r-1} \leq v_r, \label{3Eqa} \\
a_{u_r} > a_{v_r}~~ &\mbox{and} & ~~  b_{u_r} < b_{v_r}. \label{3Eqc}
\end{eqnarray}

If $p_{r-1}$, $q_{r-1}$, $u_r$ and $v_r$ are already constructed for some $ r>0$, then we continue by defining
\begin{eqnarray}
 p_r &:=& \max \{k \leq p_{r-1} ~|~ a_k> a_{v_r}\} \geq u_r, \label{3Eqp} \\
q_r &:=& \min \{k \geq q_{r-1} ~|~ b_k> b_{u_r}\} \leq v_r. \label{3Eqq}
\end{eqnarray}

Note that if $p_r<M$, we have 
\begin{equation}  \label{3Eq7}
 a_{p_r+1} \leq a_{v_r},
\end{equation}
and similarly
\begin{equation} \label{3Eq8}
b_{q_r-1} \leq  b_{u_r}
\end{equation}
if $q_r>N$.

By construction we have $i_{p_r} \prec i_{v_r}$ and $i_{u_r} \prec i_{q_r}$. According to Lemma \ref{3L4} there exists
a connecting gap for each of these pairs, and we investigate whether one of these gaps is suitable as the auxiliary element $x$.
If both of them fail to fit, then we continue our recursive procedure:

\begin{enumerate}
\item[(A)] The pair $i_{p_r} \prec i_{v_r}$: Insertion of  $x:=\alpha \beta -a_{v_r}\alpha-b_{p_r}\beta$ allows
to split the couple $[I,J]$ in the following way:

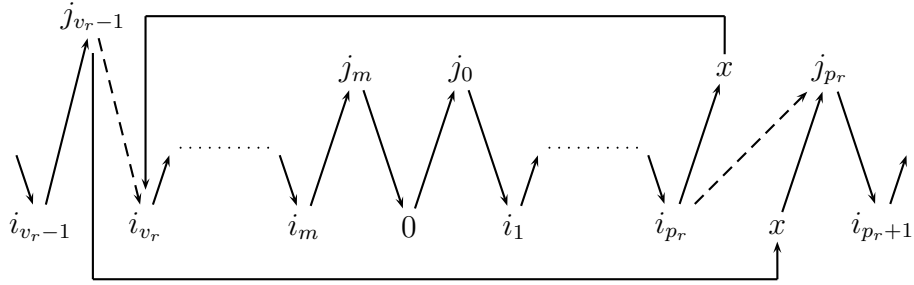
\psset{unit=0.7cm}
\begin{pspicture}(1.0,-1.5)(10,5)

\rput(0,0){\rnode{iv1}{$i_{v_{r}-1}$}}
\rput(2,0){\rnode{iv}{$i_{v_r}$}}
\pnode(2,0.7){ia}
\pnode(1,-1){ib}
\rput(5,0){\rnode{im}{$i_m$}}
\rput(7,0){\rnode{i0}{$0$}}
\rput(9,0){\rnode{i1}{$i_1$}}
\rput(12,0){\rnode{ip}{$i_{p_r}$}}
\rput(14,0){\rnode{ix}{$x$}}
\pnode(14,-1){ic}
\rput(16,0){\rnode{ip1}{$i_{p_{r}+1}$}}

\pnode(-0.5,1.5){ja}
\rput(1,4){\rnode{jv}{$j_{v_{r}-1}$}}
\pnode(1,3.3){jvv}
\pnode(2,4){jz}
\pnode(2.5,1.5){jb}
\pnode(4.5,1.5){jc}
\rput(6,3){\rnode{jm}{$j_m$}}
\rput(8,3){\rnode{j0}{$j_0$}}
\pnode(9.5,1.5){jd}
\pnode(11.5,1.5){je}
\rput(13,3){\rnode{jx}{$x$}}
\pnode(13,4){jy}
\rput(15,3){\rnode{jp}{$j_{p_r}$}}
\pnode(16.5,1.5){jf}

\ncline[nodesep=3pt]{->}{ja}{iv1}
\ncline[nodesep=3pt]{->}{iv1}{jv}
\ncline[nodesep=3pt,linestyle=dashed]{->}{jv}{iv}
\ncline[nodesep=3pt]{->}{iv}{jb}
\ncline[nodesep=3pt,linestyle=dotted]{jc}{jb}
\ncline[nodesep=3pt]{->}{jc}{im}
\ncline[nodesep=3pt]{->}{im}{jm}
\ncline[nodesep=3pt]{->}{jm}{i0}
\ncline[nodesep=3pt]{->}{i0}{j0}
\ncline[nodesep=3pt]{->}{j0}{i1}
\ncline[nodesep=3pt]{->}{i1}{jd}
\ncline[nodesep=3pt,linestyle=dotted]{jd}{je}
\ncline[nodesep=3pt]{->}{je}{ip}
\ncline[nodesep=3pt,linestyle=dashed]{->}{ip}{jp}
\ncline[nodesep=3pt]{->}{jp}{ip1}
\ncline[nodesep=3pt]{->}{ip1}{jf}

\ncline[nodesep=3pt]{->}{ip}{jx}
\ncline[nodesepA=3pt]{jx}{jy}
\ncline{jz}{jy}
\ncline{->}{jz}{ia}

\ncline[nodesep=3pt]{->}{ix}{jp}
\ncline{jvv}{ib}
\ncline{ib}{ic}
\ncline[nodesepB=3pt]{->}{ic}{ix}

\end{pspicture}

\noindent The first part
\[
  \big[(i_0, \ldots , i_{p_r},i_{v_r}, \ldots , i_m),(j_0, \ldots, j_{p_{r}-1}, x, j_{v_r}, \ldots , j_m) \big]
\]
is a fundamental couple. If its counterpart
\[
 \big[(x, i_{p_r+1}, \ldots , i_{v_r-1}),(j_{p_r}, \ldots , j_{v_r-1})\big]
\]
is a balanced couple too, then the induction hypothesis applies to it since it is of length 
\[1+ v_r-1 -(p_r+1) +1 = v_r-p_r \leq m-1.\]
 The required congruences are satisfied, hence it remains to check whether
$x \leq j_{p_r}, j_{v_r-1}$. The first inequality 
\[x \leq \alpha \beta - a_{p_{r}+1}\alpha - b_{p_r} \beta =j_{p_r}\]
 is clear for $p_r<M$ because of (\ref{3Eq7}), and if $p_r=M$ it holds since 
 in this case one has
 $j_{M} \geq 2 \alpha \beta - a_{M+1}\alpha-b_{M}\beta$, compare the proof of Lemma \ref{3L7}.
 Similarly, the second inequality is clear if $v_r=N$, since in this case we have
 $j_{v_r-1} \geq 2 \alpha \beta - a_{v_r}\alpha-b_{v_r-1}\beta$.
Otherwise, i.~e., if $v_r>N$, then $j_{v_r-1}$ is a gap. We have
\[ j_{v_r-1} = \alpha \beta -a_{v_r} \alpha -b_{v_r-1}\beta, \]
hence $x \le j_{v_r-1}$ if and only if $b_{p_r} \ge b_{v_r-1}$. 
Note that this inequality holds in particular if $p_r=u_r$ and $q_r=v_r$, since either $q_r=N$ or (\ref{3Eq8}) is valid.
\item[]
\item[(B)] The pair $i_{u_r} \prec i_{q_{r}}$: This case is mirror--imaged to (A). Insertion of
 $x:=\alpha \beta -a_{q_r}\alpha-b_{u_r}\beta$ yields a splitting of $[I,J]$ into the fundamental couple
\[
 \big[ (0,i_{u_r}, i_{q_r}, \ldots , i_{m}),(j_0, \ldots, j_{u_{r}-1}, x, j_{q_r}, \ldots , j_m) \big] 
 \]
 and its counterpart
\[
\big[ (x, i_{u_r+1}, \ldots , i_{q_r-1}),(j_{u_r}, \ldots , j_{q_r-1}) \big].
\]
Again we are done if the latter is a balanced couple. Here we have to check the inequalities $x \leq j_{q_{r-1}}, j_{u_r}$. 
The first of these 
\[x \leq \alpha \beta - a_{q_{r}}\alpha - b_{q_{r}-1} \beta \]
follows from (\ref{3Eq8}) if $q_r>N$, and in the case $q_r=N$ one has
$j_{N-1} \geq 2 \alpha \beta - a_{N}\alpha-b_{N-1}\beta$; similarly, the second inequality is clear if $v_r=M$, since then one has 
$j_{u_r} \geq 2 \alpha \beta - a_{u_r+1}\alpha-b_{u_r}\beta$. Otherwise, i.~e., if $u_r<M$, we have
\[
j_{u_r} = \alpha \beta -a_{u_r+1} \alpha -b_{u_r}\beta,
\] 
thus $x \leq j_{u_r}$ if and only if $a_{q_r} \geq a_{u_r+1}$.
\item[]
\item[(C)] By the previous discussion it remains to deal with the following situation: 
\[
b_{p_r} \overset{\mathrm{(A)}}{<} b_{v_r-1}, ~v_r>N ~~\mbox{and}~~ 
 a_{q_r} \overset{\mathrm{(B)}}{<} a_{u_r+1}, ~u_r<M. 
\]
We continue by defining the next elements of the sequences $(u_r)$ and $(v_r)$ by
\begin{eqnarray*} 
 u_{r+1} &:= &\left\{\begin{array}{cl} u_r & \mbox{~~for}~~ u_r=p_r \\ u_r +1 & \mbox{~~otherwise.}\end{array}\right. \\
& & \\
v_{r+1} &:=& \left\{\begin{array}{cl} v_r & \mbox{~~for}~~ v_r=q_r \\ v_r -1 & \mbox{~~otherwise.}\end{array}\right.
\end{eqnarray*}
\noindent Note that we cannot have $u_{r+1}=u_r$ and $v_{r+1} = v_r$ simultaneously, since the case $u_r=p_r$ and $v_r=q_r$
is covered by (A). \medskip

It is easy to see that the inequalities (\ref{3Eqb}) -- (\ref{3Eqc}) also hold for $r+1$: For the first and the second this is trivial
by our assumption resp. by definition of $(u_r)$, $(v_r)$. Since $u_{r+1} \leq u_r+1$
and $v_{r+1} \geq v_r-1$ we have $a_{u_{r+1}} \geq a_{u_r+1}$ and $b_{v_{r+1}} \geq b_{v_r-1}$; together with
(\ref{3Eqa}) this implies (\ref{3Eqc}). Hence we may continue with the construction of $p_{r+1}$ and $q_{r+1}$.
\medskip

\noindent By construction, it is clear that this recursive procedure will eventually terminate, namely with one of
the cases $u_r=M$, $u_r=p_r \,\wedge \,v_r=q_r$ or $v_r=N$, which are covered by the discussion above. \BE

\end{enumerate}  

\subsection{Proof of the main result}
After the previous rather technical subsection we return to the proof of the main result, which now finally can be completed
with the aid of Theorem \ref{Th:bl}:

Let $H(t):= \sum_{n=0}^{\infty}h_nt^n = \frac{Q(t)}{(1-t^{\alpha} ) (1-t^{\beta} )}$ be a nonnegative Laurent series 
 satisfying condition {\rm ($\star$)} with $h_0 >0$. We want to show that at least one of the series 
\[ H_{\alpha}(t) = H(t) - \frac{1}{1-t^{\alpha}}, \hspace{2cm} H_{\beta}(t) = H(t) - \frac{1}{1-t^{\beta}}\]
is nonnegative and satisfies {\rm ($\star$)} as well. Since by Proposition \ref{3P3} at least one of the numbers
\begin{eqnarray*}
 c_{\alpha} &:= &\min \{ h_{r \alpha} ~|~ r>0 \}\\
c_{\beta} &:=& \min \{ h_{r \beta} ~|~ r>0 \} 
\end{eqnarray*}
is positive, there are two cases: If only one of these series, say $H_{\beta}$, is nonnegative, then we have to show that the
 $\beta$--critical inequalities hold strictly. If both series are nonnegative, then we have to show that all critical inequalities of 
one type hold strictly. We begin with the first case.

\begin{Prop} \label{3P4}
 Let $H(t):= \sum_{n=0}^{\infty}h_nt^n$ be a nonnegative Laurent series 
 satisfying condition {\rm ($\star$)} and $h_0 >0$. If $c_{\alpha} = 0$ (resp.~$c_{\beta} = 0$), then the $\beta$--critical
(resp.~the $\alpha$--critical) inequalities hold strictly.
\end{Prop}

\dem~ Assume $c_{\alpha} = 0$, thus $h_{r \alpha}=0$ for some $0 < r < \beta$. Let
\[\sum_{i \in I} h_{n+i} \leq \sum_{j \in J} h_{n+j}\]
be a $\beta$--critical inequality, so  $n+i_p = -s' \beta$ for 
some $0 \leq p \leq m$, $s' >0$.  We define a balanced couple
\[
[I',J'] := \left[ (i_k+n)_{k=0}^m, (j_k+n)_{k=0}^m \right].
\]
Choose some integer $\ell < r \alpha, j_p'$ such that $\ell \equiv r \alpha \!\mod \beta$ and $\ell \equiv j_p' \!\mod \alpha$.
We construct another balanced couple $[I'',J'']$ by replacing $i_p'=-s' \beta$ with the sequence 
$0 \to r \alpha \to \ell$, as the following picture illustrates:

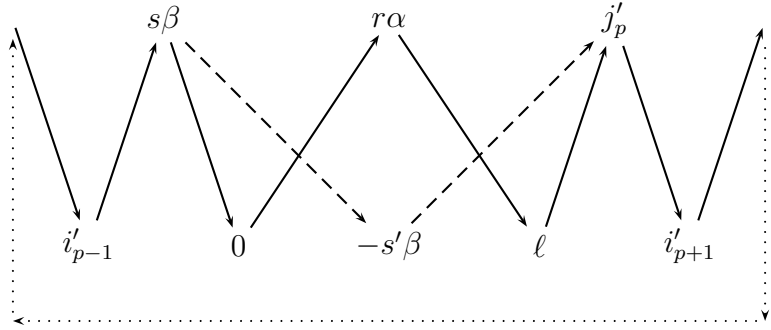
\begin{pspicture}(-2.5,-1.5)(10,3.5)

\pnode(-1,3){ia}
\pnode(-1,-1){ic}
\rput(0,0){\rnode{im}{$i_{p-1}'$}}
\rput(2,0){\rnode{n}{$0$}}
\rput(4,0){\rnode{sb}{$-s' \beta$}}
\rput(6,0){\rnode{l}{$\ell$}}
\rput(8,0){\rnode{ip}{$i_{p+1}'$}}
\pnode(9,3){ib}
\pnode(9,-1){id}

\rput(1,3){\rnode{jp1}{$s \beta$}}
\rput(4,3){\rnode{r}{$r \alpha $}}
\rput(7,3){\rnode{jp}{$j_p'$}}

\ncline[nodesep=3pt]{->}{ia}{im}
\ncline[nodesep=3pt]{->}{im}{jp1}
\ncline[nodesep=3pt]{->}{jp1}{n}
\ncline[nodesep=3pt]{->}{n}{r}
\ncline[nodesep=3pt]{->}{r}{l}
\ncline[nodesep=3pt]{->}{l}{jp}
\ncline[nodesep=3pt]{->}{jp}{ip}
\ncline[nodesep=3pt]{->}{ip}{ib}
\ncline[nodesep=3pt,linestyle=dashed]{->}{jp1}{sb}
\ncline[nodesep=3pt,linestyle=dashed]{->}{sb}{jp}
\ncline[nodesepB=7pt,linestyle=dotted]{->}{ic}{ia}
\ncline[nodesepA=7pt,linestyle=dotted]{->}{ib}{id}
\ncline[linestyle=dotted]{->}{id}{ic}
\end{pspicture}

By Theorem \ref{Th:bl} we have
$ \sum_{i \in I''} h_{i} \leq \sum_{j \in J''} h_{j}$.
This implies
\[ \sum_{i \in I} h_{n+i} < h_0 + h_{\ell} + \sum_{i \in I'} h_i =
 \sum_{i \in I''} h_i \leq \sum_{j \in J''} h_j = h_{r \alpha} + \sum_{j \in J'} h_j
= \sum_{j \in J} h_{n+j},
\]
so the original $\beta$--critical inequality holds strictly. The case $c_{\beta}=0$ is treated analogously. \BE \lz

The basic idea for solving the second case is quite similar: 

\begin{Prop} \label{3P5}
 Let $H(t):= \sum_{n=0}^{\infty}h_nt^n $ be a nonnegative Laurent series satisfying condition {\rm ($\star$)} 
and $h_0 >0$. If $c_{\alpha}, c_{\beta} \neq 0$, then the $\alpha$--critical or the $\beta$--critical inequalities hold strictly.
\end{Prop}

\dem~ Assume on the contrary that there is a non--strict $\alpha$--critical inequality, i.~e.
\[\sum_{i \in I_{\alpha}} h_{n_{\alpha}+i} = \sum_{j \in J_{\alpha}} h_{n_{\alpha}+j}\]
as well as a non--strict $\beta$--critical inequality 
\[\sum_{\tilde{\imath} \in I_{\beta}} h_{n_{\beta}+\tilde{\imath}} 
= \sum_{\tilde{\jmath} \in J_{\beta}} h_{n_{\beta}+\tilde{\jmath}} \]
where $[I_{\alpha},J_{\alpha}] = \left[ (i_k)_{k=0}^{m_{\alpha}}, (j_k )_{k=0}^{m_{\alpha}} \right]$ and 
$[I_{\beta},J_{\beta}] = \left[ (\tilde{\imath}_k)_{k=0}^{m_{\beta}}, (\tilde{\jmath}_k )_{k=0}^{m_{\beta}} \right]$. By
definition of a critical inequality there exist
 $0 \leq p \leq m_{\alpha}$, $0 \leq q \leq m_{\beta}$, and $r',s' >0$ such that
$n_{\alpha}+i_p = -r' \alpha$ and $n_{\beta}+\tilde{\imath}_q = -s' \beta$.   We define balanced couples 
\begin{eqnarray*}
 [\hat{I},\hat{J} ] &:=& \left[ (i_k+n_{\alpha})_{k=0}^{m_{\alpha}}, (j_k+n_{\alpha})_{k=0}^{m_{\alpha}} \right] \\
\left[\check{I},\check{J}\right] &:=& \left[ (\tilde{\imath}_k+n_{\beta})_{k=0}^{m_{\beta}}, 
(\tilde{\jmath}_k+n_{\beta})_{k=0}^{m_{\beta}} \right]. 
\end{eqnarray*}
We construct another balanced couple $[I,J]$ by glueing together $[\hat{I},\hat{J} ]$ and $\left[\check{I},\check{J}\right]$,
as illustrated in the following picture:

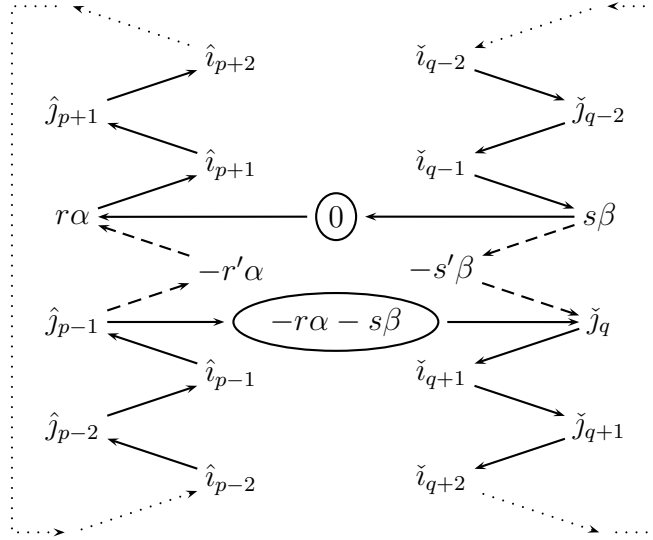
\psset{unit=0.7cm}
\begin{pspicture}(-4,-2.5)(10,9)
\pnode(0,-2){jn2}
\rput(0,0){\rnode{jn1}{$\hat{\jmath}_{p-2}$}}
\rput(0,2){\rnode{jn}{$\hat{\jmath}_{p-1}$}}
\rput(0,4){\rnode{j0}{$r \alpha$}}
\rput(0,6){\rnode{j1}{$\hat{\jmath}_{p+1}$}}
\pnode(0,8){j2}

\rput(3,-1){\rnode{in1}{$\hat{\imath}_{p-2}$}}
\rput(3,1){\rnode{in}{$\hat{\imath}_{p-1}$}}
\rput(3,3){\rnode{i0}{$-r' \alpha$}}
\rput(3,5){\rnode{i1}{$\hat{\imath}_{p+1}$}}
\rput(3,7){\rnode{i2}{$\hat{\imath}_{p+2}$}}

\ncline[nodesepB=3pt,linestyle=dotted]{->}{jn2}{in1}
\ncline[nodesep=3pt]{->}{in1}{jn1}
\ncline[nodesep=3pt]{->}{jn1}{in}
\ncline[nodesep=3pt]{->}{in}{jn}
\ncline[nodesep=3pt,linestyle=dashed]{->}{jn}{i0}
\ncline[nodesep=3pt,linestyle=dashed]{->}{i0}{j0}
\ncline[nodesep=3pt]{->}{j0}{i1}
\ncline[nodesep=3pt]{->}{i1}{j1}
\ncline[nodesep=3pt]{->}{j1}{i2}
\ncline[nodesepA=3pt,linestyle=dotted]{->}{i2}{j2}

\pnode(10,-2){jt2}
\rput(10,0){\rnode{jt1}{$\check{\jmath}_{q+1}$}}
\rput(10,2){\rnode{jt0}{$\check{\jmath}_q$}}
\rput(10,4){\rnode{jtm}{$s \beta$}}
\rput(10,6){\rnode{jtm1}{$\check{\jmath}_{q-2}$}}
\pnode(10,8){jtm2}

\rput(7,-1){\rnode{it2}{$\check{\imath}_{q+2}$}}
\rput(7,1){\rnode{it1}{$\check{\imath}_{q+1}$}}
\rput(7,3){\rnode{it0}{$-s' \beta$}}
\rput(7,5){\rnode{itm}{$\check{\imath}_{q-1}$}}
\rput(7,7){\rnode{itm1}{$\check{\imath}_{q-2}$}}

\ncline[nodesepB=3pt,linestyle=dotted]{->}{jtm2}{itm1}
\ncline[nodesep=3pt]{->}{itm1}{jtm1}
\ncline[nodesep=3pt]{->}{jtm1}{itm}
\ncline[nodesep=3pt]{->}{itm}{jtm}
\ncline[nodesep=3pt,linestyle=dashed]{->}{jtm}{it0}
\ncline[nodesep=3pt,linestyle=dashed]{->}{it0}{jt0}
\ncline[nodesep=3pt]{->}{jt0}{it1}
\ncline[nodesep=3pt]{->}{it1}{jt1}
\ncline[nodesep=3pt]{->}{jt1}{it2}
\ncline[nodesepA=3pt,linestyle=dotted]{->}{it2}{jt2}

\rput(5,4){\ovalnode{nu}{$0$}}
\rput(5,2){\ovalnode{ab}{$-r \alpha -s \beta$}}

\ncline[nodesep=3pt]{->}{jtm}{nu}
\ncline[nodesep=3pt]{->}{nu}{j0}

\ncline[nodesep=3pt]{->}{jn}{ab}
\ncline[nodesep=3pt]{->}{ab}{jt0}

\ncbar[nodesep=3pt,angle=180,arm=20pt,linestyle=dotted]{->}{j2}{jn2}
\ncbar[nodesep=3pt,arm=20pt,linestyle=dotted]{->}{jt2}{jtm2}

\end{pspicture}

By Theorem \ref{Th:bl}, we have
$ \sum_{i \in I} h_i \leq \sum_{j \in J} h_j$.
Since $h_{-r' \alpha} = h_{-s' \beta} =0$ we have
\[
 \sum_{i \in \hat{I}} h_i + \sum_{i \in \check{I}} h_i 
< h_0  + \sum_{i \in \hat{I}} h_i + \sum_{i \in \check{I}} h_i
 = \sum_{i \in I} h_i \leq \sum_{j \in J} h_j =\sum_{j \in \hat{J}} h_j + \sum_{j \in \check{J}} h_j,\]
but this contradicts
\begin{align*}
&\sum_{i \in \hat{I}} h_i + \sum_{i \in \check{I}} h_i \\
=~~&\sum_{i \in I_{\alpha}} h_{n_{\alpha}+i} + \sum_{\tilde{\imath} \in I_{\beta}} h_{n_{\beta}+\tilde{\imath}} 
=\sum_{j \in J_{\alpha}} h_{n_{\alpha}+j} +\sum_{\tilde{\jmath} \in J_{\beta}} h_{n_{\beta}+\tilde{\jmath}} \\
	=~~&\sum_{j \in \hat{J}} h_j + \sum_{j \in \check{J}} h_j. 
\end{align*} 
\BE

After these final preparatory steps we are ready to prove the essential assertion, Proposition \ref{3P2}.

{\it Proof of Prop.~\ref{3P2}:} Since $H$ is of dimension $1$, there exists an integer $N$ such that
$h_n = h_{n+ \alpha \beta}$ holds for all $n \geq N$. Then the sum  $ \sum_{k=n}^{n + \alpha \beta -1} h_k$
has the same value for every $n \geq N$; we denote this value by $\sigma(H)$. 

We prove the assertion by induction on $s:= \sigma(H)$, starting with the vacuous case $s=0$. For $s>0$ we may assume $h_0 >0$
and $h_k =0$ for $k <0$.  Let $c_{\alpha}$ and $c_{\beta}$ be defined as above. We distinguish two cases: If $c_{\alpha}$ vanishes, 
then, by Propositions \ref{3P3} and \ref{3P4}, $H_{\beta}(t)$ is a nonnegative series satisfying condition ($\star$).
 Since $\sigma(H_{\beta}) < \sigma(H)$ we are done; the same argument works with $\alpha$ and $\beta$ interchanged.
If $c_{\alpha}, c_{\beta} >0$, then both series $H_{\alpha}$, $H_{\beta}$ are nonnegative, and at least one of them also 
satisfies condition ($\star$) by Proposition \ref{3P5}, so we may apply the induction hypothesis to it. \BE \lz

As mentioned above, this result implies the converse of Proposition \ref{3P1} for any $R$--module; therefore our
main result, Theorem \ref{3T1}, is completely proven. 

The closing example of this section confirms the importance of Proposition \ref{3P5}.
\begin{Exa}
 Let $\alpha =3$ and $\beta=4$. For 
\[ H(t) := \frac{1+t+t^6+t^7+t^8}{1-t^3} = 1 + t +  0t^2 + t^3 + t^4 + 0t^5 + t^6 + \cdots , \]
we have $c_4 =1$, but not all the $4$--critical inequalities hold strictly. Hence there 
exists no decomposition of $H$ into summands $\frac{t^k}{1-t^3}$ and $\frac{t^k}{1-t^4}$ containing $\frac{1}{1-t^4}$.
\end{Exa}
\dem~ Obviously we have $h_{4r} \geq 1$ for all $r \geq 0$, but the
$4$--critical inequality
\[ h_{-4} + h_1 \leq h_4 + h_5 \]
does not hold strictly. Therefore
\[ H_4(t) := H(t) - \frac{1}{1-t^4} = 0 + t + 0t^2 + 0t^3 + 0t^4 + 0t^5 + t^6 + \cdots \]
is nonnegative, but does not satisfy condition ($\star$), and so $\nu(H_4) =0$; the latter is easily seen,
since neither $H_4(t) - \frac{t}{1-t^3}$ nor $H_4(t) - \frac{t}{1-t^4}$ is nonnegative. \\ ~ \BE

\subsection{Remarks}
i) We point out that Theorem \ref{3T1} is also valid in the degenerate case  $\alpha =1$.
Since the semigroup $\langle 1,\beta  \rangle = \NN_0$ has no gaps at all, 
condition ($\star$) collapses to the single inequality $h_n \leq h_{n+ \beta}~ \forall n \in \ZZ$.
This criterion could be deduced alternatively by applying \cite[Thm.~2.1]{ju}
to the nonnegative series $(1-t^{\beta})H_M(t)$.  \medskip

ii) The case of $\deg(X)$ and $\deg(Y)$ having a common divisor $>1$ can be reduced to the case of coprime
degrees by standard methods. Hence Theorem \ref{3T1} provides a criterion for positive Hilbert depth 
also in the general case: 

Let $\deg(X) = \alpha' = \alpha \delta$ and $\deg(Y) = \beta' =\beta \delta$ with $\delta>1$ and $\gcd(\alpha, \beta) = 1$.
Any finitely generated graded $R$--module $M = \bigoplus_n M_n$  decomposes into a direct sum of Veronese submodules
\[ M= \bigoplus_k M^{(k)}, ~~~\mbox{where}~~~ M^{(k)} := \bigoplus_{n=0}^{\infty} M_{n\delta +k},~~~ k= 0, \ldots ,\delta-1. \]
We change the grading of $R$ and $M^{(k)}$ by setting $R_{n\delta}$ resp.~$M_{n\delta +k}$ as the $n$th component in the
new grading. Then $M^{(k)}$ is still a graded $R$--module. Rewriting the conditions for positive Hilbert
depth given by Theorem \ref{3T1} in terms of the original grading yields \medskip

($\star_{k}$) \hspace*{\fill} $\sum_{i \in I} h_{id+n\delta+k} \leq \sum_{j \in J} h_{jd+n\delta+k}$  ~~ 
for all $n \in \ZZ,~ [I,J] \in \mathcal{F}_{\alpha,\beta}.$ \hspace*{\fill}
\medskip

Since
\[ 
\hdp(M) = \min_{k} \left\{ \hdp\left( M^{(k)} \right) \right\}, 
\]
we have $\hdp(M)> 0$ if and only if $H_M$ satisfies conditions ($\star_{k}$) for
$k = 0, \ldots , \delta-1$. 

iii) Theorem \ref{3T1} also holds for modules over a larger polynomial ring
 $\FF[X_1, \ldots, X_r,Y_1, \ldots ,Y_s]$ where each variable is assigned one of two coprime degrees:
The proof given above can be extended to a proof by induction on the dimension of the module, since a
reductive step similar to (\ref{3Eq4}) also works for higher dimensions. 

iv) Let $M$ be a finitely generated graded $R$--module of positive depth. As explained above, Theorem \ref{2T1}
implies that $H_M$ satisfies condition ($\star$), but from this argument it is not immediately clear why the existence of
 an $M$--regular element forces these inequalities. 
The only obvious exception is the minimal inequality $h_n \leq h_{n + \alpha \beta}$. There is also 
an alternative explanation for one special inequality with maximal number of terms: The condition
\begin{equation} \label{3EQ}
 h_n + h_{n+1} + \cdots  + h_{n+ \alpha -1} 
 \leq h_{n + \beta} + h_{n+ \beta + 1} + \cdots + h_{n+ \beta + \alpha -1}
\end{equation}
 can be deduced as follows.

\noindent Let $S = \FF[U,V]$ be the standard graded polynomial ring, then we may identify $R$ with the 
subalgebra $\FF[U^{\alpha},V^{\beta}]$ of $S$, and in this sense $S$ is a finite free $R$--module with basis
$\{ U^iV^j~|~ 0\leq i < \alpha,~ 0 \leq j < \beta\}$. Hence $\tilde{M} := M \otimes_R S$ is a finite graded $S$--module
of the same depth as $M$ with
\[ H_{\tilde{M}}(t) = \left(\sum_{i=0}^{\alpha-1}t^i \cdot \sum_{j=0}^{\beta-1}t^j \right)\cdot H_M(t) =: \sum_n \tilde{h}_nt^n.\]
Since $\dep_S(\tilde{M})>0$ we have $\pos(\tilde{M}) >0$, i.~e.~$\tilde{h}_n \leq \tilde{h}_{n+1}$ for all $n \in \ZZ$, and
rewriting this inequality in terms of $h_n$ yields exactly (\ref{3EQ}).



\begin{thebibliography}{99}
   
\bibitem{ju}  J.~Uliczka, Remarks on Hilbert Series of Graded Modules over Polynomial Rings,
              Manuscr. math.~132~(2010)~159--168.
              
\bibitem{kr} M.~Kreuzer, L.~Robbiano, Computational Commutative Algebra 2, Springer, 2005.

\bibitem{ls} L.~Smith, Polynomial Invariants of Finite Groups, second edition, A.~K.~Peters Ltd., Wellesley, MA, 1997.

\bibitem{st}  R.~P.~Stanley, Enumerative Combinatorics vol.~I, Wadsworth \& Brooks/Cole, 1986.

\bibitem{bh}  W.~Bruns, J.~Herzog, Cohen-Macaulay rings, revised edition, Cambridge University Press, Cambridge, 1998.

\bibitem{bku}  W.~Bruns, Chr.~Krattenthaler, and J.~Uliczka, Stanley decompositions and Hilbert depth in the Koszul complex,
               J.~of~Comm.~Alg.~ 2 ~(2010) ~327--357. 

\bibitem{no} W.~Bruns, B.~Ichim, Normaliz: algorithms for rational cones and affine monoids,
             J.~of~Algebra~324~(2010) ~1098--1113.
             
\bibitem{rg} J.~C.~Rosales, P.~A.~Garc\'{\i}a S\'{a}nchez, Numerical semigroups, Dev. in Math. vol. 20, Springer, Berlin-Heidelberg-New York, 2010.

\bibitem{rosales} J.~C.~Rosales, Fundamental gaps of numerical semigroups generated by two elements,
                  Linear Algebra and its Appl.~405~(2005)~200--208.

\bibitem{si} S.~Ilhan, On Ap\'ery Sets of Symmetric Numerical Semigroups,
             Int.~Math.~Forum~1~(2006) ~481--484. 

\end{thebibliography}
\end{document}